\documentclass[12pt,leqno]{article}
\usepackage{amsmath,amsthm,amssymb}

\def\init{\setcounter{equation}{0}}
\newtheorem{theorem}{Theorem}[section]
\newcommand{\R}{\mathbb{R}}

\newcommand{\Z}{\mathbb{Z}}

\newcommand{\e}{{\varepsilon}}

\newcommand{\rw}{\rightarrow}

\usepackage{tikz}
\usetikzlibrary{decorations.pathreplacing}

\usepackage{verbatim}

\title{Determination  of black holes by boundary measurements}

\author{G.Eskin, \ \ \  Department of Mathematics, UCLA,\\ Los Angeles,
CA 90095-1555, USA. \ E-mail: eskin@math.ucla.edu
}

\begin{document}

\maketitle

\begin{abstract}
For  a wave  equation  with time-independent   Lorentzian   metric  consider   an initial-boundary value
problem   in $\R\times \Omega$,  where  $x_0\in  \R$  is  the  time  variable  and  $\Omega$  is  a bounded domain 
in $\R^n$.   Let  $\Gamma\subset\partial\Omega$  be a  subdomain  of  $\partial\Omega$.
We say  that  the boundary  measurements  are   given  on  $\R\times\Gamma$  if we  know   the  Dirichlet  and  Neumann data  
on  $\R\times \Gamma$.   
The  inverse  boundary  value  problem  consists  of  recovery   of the metric  from the boundary  data.   In author's 
previous  works  a localized  variant  of the boundary  control method was  developed  that allows  the  recovery of  the metric  locally
in a neighborhood  of any point  of $\Omega$   where  the  spatial  part  of the wave operator  is elliptic.  
 This allow the recovery  of the
metric in the exterior  of the  ergoregion.

Our goal  is to recover  the black hole.  In  some cases  the ergoregion  coincides  with  the black  hole.
 In the case of two space dimensions  we recover  the black  hole   inside  the ergoregion  
assuming  that the ergosphere,  i.e.  the boundary  of the  ergoregion,  is not  characteristic  at any  point  
of the ergosphere.
 
\end{abstract}

\section{Introduction}
\init

Let  $(x_0,x_1,x_2,...,x_n)\in  \R\times\R^n, \ x_0\in \R$   be  the time variable,   $x=(x_1,x_2,...,x_n)\in \R^n$. 
Consider  the Lorentzian metric  in  $\R\times\R^n$
\begin{equation}															\label{eq:1.1}
\sum_{j,k=0}^n g_{jk}(x)dx_jdx_k
\end{equation}
with the signature $(+1,-1,...,-1)$.    Let  $g(x)=\det [g_{jk}(x)]_{j,k=0}^n$   and   let $  [g^{jk}(x)]_{j,k=0}^n$   be the inverse  to the  metric tensor
$ [g_{jk}(x)]_{j,k=0}^n$.   We assume  that  the metric  does not depend   on $x_0\in \R$.

Let
\begin{equation}															\label{eq:1.2}
Lu =\sum_{j,k=0}^n\frac{1}{\sqrt{(-1)^ng(x)}}\frac{\partial}{\partial x_j}\Big(\sqrt{(-1)^ng(x)}g^{jk}(x)\frac{\partial}{\partial x_k}u(x)\Big)=0
\end{equation}
be  the wave  equation  corresponding  to the metric  (\ref{eq:1.1}).

We assume that  
\begin{equation}															\label{eq:1.3}
g^{00}(x)>0.
\end{equation}
Let  $\Omega$  be a bounded domain  in  $\R^n$.

We shall consider  the following  initial-boundary value  problem in  $\R\times \Omega$  for  the equation  (\ref{eq:1.2})  
\begin{equation}																\label{eq:1.4}
u=0  \ \ \ \mbox{for}\ \ \ x_0 \ll 0,\ x\in \Omega,
\end{equation}
\begin{equation}																\label{eq:1.5} 
u\big|_{\R\times \partial\Omega}=f,
\end{equation}
where  $\R\times \partial\Omega$  is  not a characteristic surface  and $f$  has a compact  support  in  $\R\times\partial\Omega$.  
Since  $\R\times \partial\Omega$  is  not characteristic  there exists a unique  solution  of the problem 
(\ref{eq:1.2}),  (\ref{eq:1.4}),  (\ref{eq:1.5})
(cf.  [14]).

Define  the Dirichlet-to-Neumann  operator $\Lambda f$  as
\begin{equation}																\label{eq:1.6}
\Lambda f=\sum_{j,k=0}^n g^{jk}(x)\frac{\partial u}{\partial x_j} \nu_k(x)  \Big(\sum_{p,r=0}^n g^{pr}(x)\nu_p(x)\nu_r(x)\Big)^{-\frac{1}{2}}  
 \Big|_{\R\times\partial\Omega},
 \end{equation}
  $u(x_0,x)$  is  the solution  of (\ref{eq:1.2}),  $\nu_0=0,(\nu_1,...,\nu_n)$  is  the  unit  outward  normal  to  $\partial \Omega$.
  
  Let  $\Gamma_0$  be a subdomain  of  $\partial\Omega$.  The  inverse  boundary value  problem  on  $\R\times\Gamma_0$
  consists  of  determining  the  metric  (\ref{eq:1.1})  knowing  the  $\Lambda f$  on  $\R\times\Gamma_0$
  for  all $f$   with the compact support in  $\R\times\overline\Gamma_0$.
  
  A  powerful boundary  control method  for  solving  hyperbolic  inverse  problems was discovered  by  M.Belishev  (cf.  [2]),  and  further
  developed in  [3],  [4],  [5],  [15],  [16],  [17]).   In author's  works  [6],  [7],  [8],  a localized  variant  of the  boundary  control
    method  was  developed  that allows  to  recover  the  metric  in a neighborhood  of any  point  of $\Omega$  where   the spatial  part
    of the wave operator  is elliptic.  These  results  are the  basis  of the  present  paper.
  
  Let  $y=\varphi(x)$  be  a diffeomorphism  of $\overline\Omega$  on  some  bounded smooth domain
  $\overline\Omega_0\subset\R^n$   and let   $a(x)\in C^\infty(\Omega),\ a(x)=0$  on $\overline\Gamma_0$.
  
  Consider the map
  \begin{equation}																\label{eq:1.7}
  (y_0,y)=\Phi(x_0,x)=(\varphi(x),x_0+a(x))
  \end{equation}
   of 
  $\overline\Omega\times\R$ onto  $\overline\Omega_0\times\R$
  such that  $\varphi(x)=x$  and  $a(x)=0$     on $\Gamma_0$.
  
  Note that	  change of variables  $y=\varphi(x),  y_0=x_0+a(x)$  does  not  change  the  DN  operator $\Lambda$.
  
  The domain  $\Delta\subset\Omega$  is called  the ergoregion  if  (see,  for example,  [22]),
  \begin{equation}															\label{eq:1.8}
  g_{00}(x)\leq 0\ \ \mbox{on}\ \ \ \Delta.
  \end{equation}
  We assume  that  $g_{00}(x)>0$  in the  exterior  of  $\Delta$.  Let
  \begin{equation}															\label{eq:1.9}
  \Delta(x)=\det[g^{jk}(x)]_{j,k=1}^n.
  \end{equation}
  Then   $g_{00}(x) =g^{-1}(x)\Delta(x)$.  Thus  (\ref{eq:1.8})  is equivalent  to the inequality
  \begin{equation}															\label{eq:1.10}
\Delta(x)\leq 0.
\end{equation}
We assume  that  $\Delta(x)=0$  is  a smooth  surface  in  $\R^n,\ \frac{\partial\Delta(x)}{\partial x}=
\big(\frac{\partial \Delta(x)}{\partial x_1},...,\frac{\partial \Delta(x)}{\partial x_n}\big)\neq 0$  when  $\Delta(x)=0$.

Now  we shall  define the black hole.

Let  $S(x)=0$  be a closed surface  in  $\R^n$  and   $\Omega_{int}$  be the  interior of  the  surface  $S(x)=0$.   We call  
the region  $\Omega_{int}$  a black hole  if no signal  (disturbance)   inside   $S(x)=0$  can reach  the exterior  of  $S(x)=0$.
Let  $S(x)=0$  be a characteristic  surface  for the equation  (\ref{eq:1.2}),
i.e.
\begin{equation}															\label{eq:1.11}  
\sum_{j,k=0}^n g^{jk}(x) S_{x_j}(x)S_{x_k}(x)=0  \ \ \mbox{when}\ \ S(x)=0.
\end{equation}
We  have  (cf.  [9])  that  $S(x)=0$
in a boundary  of a black  hole  if  $S(x)=0$  is a characteristic  surface  and
\begin{equation}															\label{eq:1.12}
\sum_{j=1}^n g^{j0}(x) S_{x_j}(x)<0 \ \ \mbox{when}\ \ S(x)=0.
\end{equation}
The  boundary  $S(x)=0$  of the  black  hole   is called  the black  hole  event  horizon.

When  the metric  (\ref{eq:1.1})   is not  a solution  of the  Einstein  equation,   i.e.    the metric  (\ref{eq:1.1})  is not related  to the 
general  theory  of relativity  (cf. [24]),    the  black  hole  $\Omega_{int}$  is called  an analogue  black hole.
In physical applications  the analogue black holes  appear  when  one  studies  the propagation of waves in a moving medium  
(cf.  [1], [20]).
An example  of the analogue  metric is the following  acoustic metric (cf.  [21],  [22]):

Consider  a fluid  flow  in a vortex  with  the  velocity  field
\begin{equation}															\label{eq:1.13}
v=(v^1,v^2)=\frac{A}{r}\hat r +\frac{B}{r}\hat\theta,
\end{equation}
where      $r=|x|,\ \hat r=\big(\frac{x_1}{|x|},\frac{x_2}{|x|}\big),\ \hat\theta=\big(-\frac{x_2}{|x|},\frac{x_1}{|x|}\big),\ A$  and $B$   are constants,
$A<0$.  When  $B\neq 0$   (\ref{eq:1.13})   is a rotating  flow.

The inverse  metric  tensor  $[g^{jk}]_{j,k=1}^2 $  has the form   
\begin{align}														\label{eq:1.14}
&g^{00}=\frac{1}{\rho c},\ g^{0j}=g^{j0}=\frac{1}{\rho c}v^j,\ 1\leq j\leq 2,
\\
\nonumber
&g^{jk}=\frac{1}{\rho c}(-c^2\delta_{ij}+v^Jv^k), \ 1\leq j,k\leq 2.
\end{align}
Here  $c$  is the  sound  speed,  $\rho$   is the density.
We shall assume,  for the simplicity,   that  $\rho=1,\ c=1$.
It was shown  in [22]   that  $\{r\leq \sqrt{A^2+B^2}\}$  is the  ergoregion  and  $\{r<|A|\}$  is the black  hole.

We conclude the introduction  by a brief  description  of the content  of this paper.

In \S2  we show  that the boundary  data given  on any  subdomain  of $\R\times\partial\Omega$  allow  to recover  the  metric  (\ref{eq:1.1})  
outside  the ergoregion.

We specify the particular case of the inverse boundary value problem when the  ergoregion  coincides with the black hole.

We  consider another  example  of analogue metric,  the Gordon  metric  ([18], [19]),  that arise when  one studies  the propagation  of light
in a moving dielectric  medium:  

Let $w=(w_1(x),w_2(x),w_3(x))$  be  the velocity  of the flow.  The Gordon metric  has the form
\begin{equation}														\label{eq:1.15}
\sum_{j,k=0}^3g_{jk}(x)dx_jdx_k,
\end{equation}
where
$g_{jk}(x)=\eta_{jk}+\big(n^{-2}(x)-1\big)v_j  v_k,\ n(x)$  is the  index  of the refraction, 
\begin{equation}													\label{eq:1.16}
v_0=\big(1-\frac{|w|^2}{c^2}\big)^{-\frac{1}{2}},\ 
v_j(x)=\big(1-\frac{|w|^2}{c^2}\big)^{-\frac{1}{2}}\frac{w_j(x)}{c},\ 1\leq j\leq 3,
\end{equation}
$\eta_{jk}$  is the Lorentz  metric.

In \S3  for $n=2$  we describe  the black hole  inside the ergosphere.

In  \S4  we  show  that knowing the metric on  the ergosphere  one can recover  the black hole inside.

In \S5  we consider  the axisymmetric    case in three  space  dimensions.

\section{Recovery of the  ergosphere  from  the boundary  measurements}
\init

Let  $\Gamma'$  be  any small subset  of  $\partial  \Omega$  and  $P_0\in  \Gamma'$.

It was proven in [6],  Theorem 2.3,  (see also  [7],  [8]),    that  knowing boundary  data  on  $[0,+\infty)\times\Gamma'$
one  can  recover,  modulo  change  of variables  (\ref{eq:1.7}),  the metric  (\ref{eq:1.1})  in $[0.+\infty)\times  V(P_0)$
where  $V(P_0)$  is a neighborhood  of  $P_0$  in  $\overline\Omega$.
The key condition for the the  validity  of  Theorem 2.3  (condition  2.3 in  [6])  is  that the spatial  part  of the equation  (\ref{eq:1.2}) 
is elliptic in $V(P_0)$,   i.e.   $V(P_0)$  is  outside  the ergosphere.  

Next  taking $P_1\in  V(P_0)$  one can recover  (\ref{eq:1.1})
in  $[0,\infty)\times V(P_1)$
where  $V(P_1)\subset \overline\Omega$  is a neighborhood   of  $P_1$.   Repeating  this  argument  infinitely  many times  
we can recover  the metric outside the ergosphere
$\Delta(x)=0$  when  $\Delta(x)$  is  the same  as in  (\ref{eq:1.9}).   Taking   the limit we can  recover the metric  on  $\Delta(x)=0$   too.
Thus we have  the following theorem:
\begin{theorem}															\label{theo:2.1}
Knowing the DN  operator  (\ref{eq:1.6})  on  $\R\times\Gamma_0$  we can  determinate  the ergosphere $\Delta(x)=0$,  and the metric  (\ref{eq:1.1}) 
on  $\Delta(x)=0$.
\end{theorem}

Note  that  the matrix  $[g^{jk}(x)]_{j,k=1}^n$  has rank $n-1$  when \linebreak $\Delta(x)=\det[g^{jk}(x)]_{j,k=1}^n=0$.    Therefore   there exists  a 
null-vector  \linebreak $e(x)=(e_1(x),...,e_n(x))$   of the matrix  $[g^{jk}(x)]_{j,k=1}^n$   smoothly  dependent  on  $x\in  \Delta(x)$
\begin{equation}																\label{eq:2.1}
\sum_{k=1}^ng^{jk}(x)e_k(x)\equiv 0,\ \ 1\leq j\leq n, \  \Delta(x)=0.
\end{equation}

We  will consider  separately  two cases:

a)   $e(x)$  is  orthogonal  to the surface  $\Delta(x)=0$  for all  $x$,  i.e.  $\Delta(x)=0$  is  characteristic  at any  $x\in\Delta$.

b)  $e(x)$  is not  orthogonal  to the surface  $\Delta(x)=0$   for any $x$,  i.e. is not characteristic  for  any $x\in \Delta$.

Note that there are many   other   cases when  $e(x)$  is  orthogonal  to  $\Delta(x)=0$  only  on a part of  $\Delta(x)=0$  but  we will
not  consider  them.

In the case  a)   we have  that  $e(x)$  is collinear  to the gradient  $\frac{\partial \Delta(x)}{\partial  x},\Delta(x)=0,$ 
for  all $x$.         Therefore  it follows  from  (\ref{eq:2.1})  that
\begin{equation}																\label{eq:2.2}
\sum_{k=1}^ng^{jk}(x)\Delta_{x_k}(x)=0, \ 1\leq j\leq n, \ \Delta(x)=0.
\end{equation}
Multiplying  (\ref{eq:2.2})  by  $\Delta_{x_j}$  and summing in $j$  we get
\begin{equation}																\label{eq:2.3}
\sum_{j,k=1}^ng^{jk}(x)\Delta_{x_j}\Delta_{x_k}=0  \ \ \mbox{when}\ \ \ \Delta(x)=0,
\end{equation}
i.e. $\Delta(x)=0$  is a characteristic  surface.  Therefore  $\Delta(x)=0$  is a boundary  either  of black hole  or a white hole
(cf. [9]).

It follows  from [9]  that  $\Delta(x)=0$  is the boundary  of a black hole  if
\begin{equation}																\label{eq:2.4}
\sum_{j=1}^n g^{0j}(x)\frac{\partial\Delta}{\partial x_j}<0\ \ \ \mbox{when}\ \ \ \Delta(x)=0.
\end{equation}
Consider,  for example,  a Schwartzchield  metric.  It has  the following form  in Cartesian  coordinates  (cf.  [23]):
\begin{equation}																\label{eq:2.5}
ds^2 = \Big(1-\frac{2m}{R}\Big)dt^2-dx^2-dy^2-dz^2-\frac{4m}{R}dtdR- \frac{2m}{R}(dR)^2,
\end{equation}
where  $R=\sqrt{x^2+y^2+z^2}$.

It is easy  to see  (cf.  [23]) that
\begin{equation}																\label{eq:2.6}
1-\frac{2m}{R}=0
\end{equation}
is simultaneously  an  ergosphere  and  a characteristic  surface.
Thus  $\{1-\frac{2m}{R}<0\}$  is  an  ergoregion  and  a black  hole.

We  shall  call  the  metric such  that  the  ergoregion  coincide  with the black hole  the  Schwartzschield type  metric.

Another example  of a Schwartzschield  metric  is  an acoustic  metric  (\ref{eq:1.13})  when  $B=0,\  A<0$.

Therefore  the problem  of recovery of the black hole by  the  boundary  measurements  is solved  for  the Schwartzshield type
metrics  since  it consists  of the recovery  of the  ergosphere.  

Further  example  of  Schwartzschield type  black  hole    appears  in   the study  of the Gordon  equation  
corresponding  to the  metric  (\ref{eq:1.16}).

The  condition  $\Delta(x)=0$  is equivalent  to  $g_{00}=0$.

Using  (\ref{eq:1.16})  we get
\begin{equation}																\label{eq:2.7}
g_{00}=1+(n^{-2}-1)v_0^2=1+\frac{(n^{-2}-1)}{1-\frac{|w|^2}{c^2}}=0.
\end{equation}
Therefore
\begin{equation}																\label{eq:2.8}
|w(x)|^2=\frac{c^2}{n^2(x)}
\end{equation}
is the ergosphere  for  the  Gordon equation. 

In the case  of the  Gordon  metric  we  have,  from  (\ref{eq:1.16})  and  (\ref{eq:2.3}):

\begin{equation}																\label{eq:2.9}
|\Delta_x|^2=\frac{n^2-1}{c^2\big(1-\frac{|w|^2}{c^2}\big)}\Big(\sum_{j=1}^nw_j(x)\Delta_{x_j}\Big)^2.
\end{equation}
Since    $|w(x)|^2=\frac{c^2}{n^2(x)}$  we get  $|\Delta_x\cdot w|^2=|\Delta_x|^2|w|^2$.
Therefore
\begin{equation}														\label{eq:2.10}
\Delta_x(x)=\alpha(x) w(x)\ \ \ \mbox{when}\ \ \ \Delta(x)=0.
\end{equation}
Since  $\Delta_x$  is  outward  normal  we have  that
$\alpha(x)<0$ when   $w(x)$  is  pointed  inside  $\Delta(x)=0$.

The  condition  (\ref{eq:2.4})  has  the  following  form  for the  Gordon  metric
\begin{equation}														\label{eq:2.11}
\sum_{j=1}^n\frac{(n^2-1)}{1-\frac{|w|^2}{c^2}}\frac{w_j}{c}\Delta_{x_j}<0.
\end{equation}
It follows  from  (\ref{eq:2.10})  that  (\ref{eq:2.11})  holds  when  $w(x)$  is  pointed  inside  $\Delta(x)=0$,
and  hence  $\Delta(x)=0$  is a boundary  of a black hole.

The case  of not  Schwartzschield  type metrics  is more  difficult.  We shall study  only  the case  of two dimensions and a metrics  
such  that any point of the ergosphere  is non-characteristic  (i.e. the case  b) ).

First  we describe  the black hole  for  this case.

\section{Description of the black hole  inside  the ergosphere  in the case of  two space dimensions}
\init
Let  $\Delta(x)=0, n=2$,  be  the ergosphere.  Assume  that the normal  to that ergosphere is not characteristic  for any  
$x\in  \{\Delta(x)=0\}$
\begin{equation}														\label{eq:3.1}  
\sum_{j,k=1}^2g^{jk}(x)\nu_j\nu_k\neq 0\ \ \ \mbox{for all}\ \ x\in \{\Delta(x)=0\},
\end{equation}
where  $(\nu_1,\nu_2)$  is the unit  normal  to  $\Delta=\{\Delta(x)=0\}$.
We assume  that  the  ergosphere  $\Delta(x)=0$  is smooth,  
i.e. $\frac{\partial \Delta}{\partial x}=\big(\frac{\partial\Delta}{\partial x_1},\frac{\partial\Delta}{\partial x_2}\big)\neq 0$
when   $\Delta(x)=0$.

As in  [11]  introduce  coordinates   $(\rho,\theta)$  where  $\rho=0$  is  the equation  of  $\Delta(x)=0,\rho=-\Delta(x)$  near  $\rho=0$.
For  the  convenience we extend  $\theta\in [0,2\pi]$  to  $\theta=\R/2\pi \Z$.  We have $0\leq\rho\leq\rho_0(\theta)$  where
$\rho=\rho_0(\theta)$  is  the  black  hole  event  horizon,  $\theta\in\R/2\pi \Z$.

Since the set $\{\Delta(x)=\det[g^{jk}(x)]_{j,k=1}^2<0\}$ is  inside  the  ergosphere,  there  are  two  characteristics  $S^\pm(x)$  such that
$$
\sum_{j,k=1}^2 g^{jk}(x)S_{x_j}^\pm S_{x_k}^\pm=0,\ \ \Delta(x)<0,
$$
or,  in $(\rho,\theta)$  coordinates,
\begin{equation}																\label{eq:3.2}
\hat g^{\rho\rho}\big(\hat S_\rho^\pm\big)^2 +2\hat g^{\rho\theta}\hat S_\rho^\pm \hat S_\theta^\pm
 +\hat g^{\theta\theta}\big(\hat S_\theta^\pm\big)^2=0,
\end{equation}
where  $\left[\begin{matrix}
\hat g^{\rho\rho}   &  \hat g^{\rho\theta}\\
\hat g^{\rho\theta} &  \hat g^{\theta\theta}
\end{matrix} \right]$
is   the matrix  $\left[\begin{matrix}
 g^{11}   &   g^{12}\\
 g^{21} &   g^{22}
\end{matrix} \right]$
in  $(\rho,\theta)$   coordinates.

We  assume that  $\hat S^\pm(\rho,\theta)$  satisfy  the following  boundary conditions
\begin{equation}															\label{eq:3.3}
\hat S^\pm(0,\theta)=\theta  \ \ \ \mbox{for any}\ \ \ \theta\in \R/2\pi\Z.
\end{equation}
\
Solving  the  quadratic  equation  (\ref{eq:3.2})  we get
\begin{equation}															\label{eq:3.4}
\hat S_\rho^\pm(\rho,\theta)=\frac{-\hat g^{\rho\theta}\pm \sqrt{-\tilde\Delta}}{\hat g^{\rho\rho}}\hat S_\theta(\rho,\theta),
\end{equation}
where  $\tilde \Delta =\hat g^{\rho\rho}\hat g^{\theta\theta} -(\hat g^{\rho\theta})^2$.

Consider the equations   for  the  null-bicharacteristics           (null  geodesics).   When  we use  the time variable  as a parameter  we have 
(cf. [9])  two families of null-geodesics:\linebreak
the  $(+)$   null-geodesics  
$
(\rho^+(x_0),\theta^+(x_0)),\ \rho^+(0)=0,\ \theta^+(0)=\theta_0,\ x_0 \geq 0,$  and   
\linebreak
the $(-)$ null-geodesics
$
(\rho^-(x_0),\theta^-(x_0)),\ \rho^-(0)=0,\ \theta^-(0)=\theta_0,\ x_0 \leq 0.$ 

Note  that  condition  (\ref{eq:3.1})  is equivalent  that  both  $(+)$  and  $(-)$  null-geodesics  are not  tangent to $\Delta(x)=0$.

Assume that the ergoregion  $\{\Delta(x)\leq 0\}$ contains  a trapped region  $O_\e$,  i.e.  a region  that both  $(+)$  and $(-)$  null-geodesics  reach   
when  $x_0\rw +\infty$  and stay there.  One of the examples  of  trapped  region is a neighborhood  of a singularity  of  the metric described  
in  [11]  (cf. (1.3),  (1.4)  in [11]).  

{\bf Remark 3.1}.  The trapped region described  above  leads  to the  existence  of the  black  hole.  
If we change the definition  of  the trapped  region  requiring that  all  $(+)$  and  $(-)$  null-geodesics  tends  to  $O_\e$  when  $x_0\rw -\infty$  we  get
that  there exists  a  white  black  hole  (cf.  [9]). 

When  $x_0\rw +\infty\ \  (\rho^+(x_0),\theta^+(x_0))$  reaches  the trapped  region  $O_\e$  and remains there  for all  large  $x_0$.

The  null-geodesics   $\gamma_-=(\rho^-(x_0),\theta^-(x_0)),\ x_0<0$,   ends  at  $((0,\theta_0)$   when      $x_0=0$  and  
$(\rho^-(x_0),\theta^-(x_0))$  cannot reach  $O_\e$  when  $x_0\rw -\infty$.   Thus    $(\rho^-(x_0),\theta^-(x_0))$  has  no limit  points  in  
$\{\Delta(x)\leq 0\}\setminus O_\e$.  

Therefore   the limiting  set  of  the  trajectories  $\{\rho^-(x),\theta^-(x_0),x_0<0\}$   is inside  $\Delta\setminus O_\e$.  
By  the Poincare-Bendixson theorem
(cf.  [13])  there  exists a  limit  cycle,  i.e.  closed  periodic solution $\gamma_0=\{ (\rho_0^-(x_0),\theta_0^-(x_0))\}$
that is  a black  hole  event  horizon.   The solution  $\gamma_-$  approaches  $\gamma_0$  spiraling  around  $\gamma_0$.  All other  
$(-)$ null-geodesics  also approach  $\gamma_0$  spiraling.

Denote  by  $\Pi$  the infinite  strip
\begin{equation}																\label{eq:3.5}
\Pi=\{0\leq \rho<\rho_0(\theta),\ -\infty<\theta<+\infty\},
\end{equation}
where  $\{\rho=\rho_0(\theta),\theta \in \R/ 2\pi\Z \}$  is the equation  of the black  hole event  horizon  $\gamma_0$.
Therefore  $\Pi$  is  traced  by all  $(-)$  null-geodesics,  and
  $\gamma_0$  is the boundary  of  $\Pi$: $\gamma_0\subset\partial\Pi$.   Thus the  following  theorem  holds:
  \begin{theorem}																\label{theo:3.1}
  Let  (\ref{eq:3.1})   holds,  i.e.  the ergosphere  $\Delta(x)=0$   is  not characteristic for all $x\in  \{\Delta(x)=0\}$.
  Suppose  the  ergoregion  $\{\Delta(x)\leq 0\}$  contains  a trapped region  $O_\e$.   Then  there exists  a black  hole
  $\overline\Pi=\{0\leq \rho\leq\rho_0(\theta),\ -\infty<\theta<+\infty\}$  and  $\{\rho=\rho_0(\theta),\theta\in  \R/2\pi\Z\}$
  is the black hole event horizon.
  \end{theorem}

In the examples  to verify  that  the  requirement  for
 the existence  of the  black  hole  inside the ergosphere  are satisfied, 
 one  need to check  that $\Delta$  is  not characteristic  when  $\Delta(x)=0$,  i.e.
 \begin{equation}																\label{eq:3.6}
 \sum_{j,k=1}^n g^{jk}(x)\frac{\partial \Delta}{\partial x_j}\frac{\partial \Delta}{\partial x_k} \neq 0  \ \ \ \mbox{when}\ \ \ \Delta(x)=0.
 \end{equation}

 In the case of the Gordon  metric  we have  that  $\Delta(x)=0$  is  equivalent  to  (cf. (\ref{eq:2.8}))
 \begin{equation}															\label{eq:3.7}
 |w(x)|^2=\frac{c^2}{n^2(x)}.
 \end{equation}

  Equation  (\ref{eq:3.6}))  in the case  of Gordon metric gives 
  \begin{align}															\label{eq:3.8}
  -|\Delta_x|^2+\Big(\frac{|w|^2}{c^2}\Big)^{-1}\big(\sum_{i=1}^2\Delta_{x_j}\frac{w_j}{c}\Big)^2\neq 0,
  \end{align}
  where we used again  (\ref{eq:3.7}).  Thus,  if  (\ref{eq:3.8})  holds,   there exists  a  black hole  inside  the  ergosphere. 
  
  Now we shall  verify  the  conditions for the existence  of the black hole  in the case  of the acoustic metric.
   
Introduce   polar coordinates   $(r,\varphi)$.   Then  the  Hamiltonian in polar  coordinates  has  the form
\begin{equation}															\label{eq:3.9}
H=\Big(\xi_0+\frac{A}{r}\xi_r+\frac{B}{r}\frac{\xi_\varphi}{r}\Big)^2-\Big(\xi_r^2+\frac{1}{r^2}\xi_\varphi^2\Big).
\end{equation}
Therefore   
  $g^{rr}=\Big(\frac{A^2}{r^2}-1\Big),\ g^{r\varphi}=\frac{AB}{r^3},\  g^{\varphi\varphi}=\frac{B^2}{r^4}-\frac{1}{r^2}$.
  
  The  ergosphere  is 
  \begin{equation}															\label{eq:3.10}
  \Delta_0=\Big(\frac{A^2}{r^2}-1\Big)\Big(\frac{B}{r^4}-\frac{1}{r^2}\Big)-\frac{A^2B^2}{r^6}=
  -\frac{A^2}{r^4}-\frac{B^2}{r^4}+\frac{1}{r^2}=0.
  \end{equation}

It is  easy to check  that
$\Delta_0$  is not  characteristic.   Therefore  a black  hole exists  inside the ergosphere.

One can show  (cf.  [9])   that  $\{r<|A|\}$  when  $A<0$  is  a black hole.

{\bf  Remark  3.2}.   Consider the case  $B=0$,   i.e.  non-rotating  fluid.  Then  the ergosphere  is   $\Delta_0(x)=\frac{1}{r^2}
\big(1-\frac{A^2}{r^2}\big)=0$.  
It was shown   in  [9]  that  $\{r<|A|\}$  is a black hole.   Thus we  have   another  example  of the Schwartzschield  type black hole. 

\section{Recovery  of  the  black  hole knowing the boundary  data on  the  ergosphere}
\init

Consider  the  characteristic  equations  (\ref{eq:3.2}).
Note that
\begin{equation}																\label{eq:4.1}
\tilde\Delta(\rho,\theta)=g^{\rho\rho}g^{\varphi\varphi}-(g^{\rho\varphi})^2\equiv   -C_1\rho,\ C_1>0.
\end{equation}
 
Let $\rho=\rho^\pm(x_0),\  \theta=\theta^\pm(x),\   x_0\geq 0,$   be  the null-geodesics  (cf.  \S3).\linebreak
Then
\begin{equation}																\label{eq:4.2}
S^\pm(\rho^\pm(x_0),\theta^\pm(x_0))\equiv S^\pm(\rho(0),\  \theta(0)),\rho(0)=0,\  \theta(0)=\theta_0.
\end{equation}
We have  that   $(\rho^+(x_0),\theta^+(x_0))$  crosses  the black hole  horizon  $\gamma_0=\{\rho=\rho_0(\theta),\theta\in \R)$   
at  some  point  $x=x_0$
and 
remains   inside  the black  hole.  The  null-geodesics   $(\rho^-(x_0),\theta^-(x_0))$  approache  the  black  hole
event  horizon  $\gamma_0$  when   $x_0\rw -\infty$.   Thus  $\gamma_0$  is  the limit  set  of  $(\rho^-(x_0),\theta^-(x_0)),x_0<0$.

Periodically  extending  $\theta\in [0,2\pi]$  to  $\theta\in (-\infty,+\infty)$  we have  that \linebreak the  $(-)$  null-geodesics  cover the strip  
$\Pi=\{0\leq \rho <\rho_0(\theta),-\infty<\theta<+\infty\}$   when $0\leq  \rho(x_0) < \rho_0,\ \theta(x_0)=\theta,  \ -\infty <\theta<+\infty$.

Let  
\begin{equation}																\label{eq:4.3}
  \sigma=S^+(\rho,\theta),\ \ \tau  =S^-(\rho,\theta)
 \end{equation}
  when $ (\rho,\theta)\in \Pi,\ \Pi=\{0\leq\rho<\rho_0(\theta), \theta\in \R^1\},$ 
$\rho=\rho_0(\theta)$  is  the black hole  event  horizon.  

Note that $Lu=0$   has  
the form  
\begin{equation}															\label{eq:4.4}
\frac{\partial^2 u}{\partial\sigma \partial\tau}=0
\end{equation}
in coordinates  $(\sigma,\tau)$.

Make a new change  of variables  
\begin{equation}														\label{eq:4.5}
y_1=\frac{\sigma+\tau}{2}=\frac{S^+(\rho,\theta)+ S^-(\rho,\theta)}{2},
\end{equation}
\begin{equation}														\label{eq:4.6}
y_2=\frac{\sigma-\tau}{2},\ \ y_1\big|_{\rho=0}=\theta,\ \ y_2\big|_{\rho=0}=0.
\end{equation}

We have that  $(y_1,y_2)=\Phi(\sigma,\tau)$  is the map  of  $\Pi$  onto  the half-plane 
$\overline\R_+^2=\{(y_1,y_2),y_1\geq 0,y_2\in (-\infty,+\infty)\}$.   Note that  $\Phi$  is  a homeomorphism  of   $\Pi$  onto
$\overline\R_+^2$,
 $\Phi$    is  the identity on  $\rho=0$.

The  closure  $\overline\Pi=\{0\leq\rho\leq \rho_0(\theta),\theta\in\R^1\}$ has the form   $\overline\Pi=\Pi\cup\gamma_0$. 
Thus  the black  hole event  horizon $\gamma_0$  belongs  to  the closure  of $\Pi$.

Suppose we  have another  metric  $g_1$  in $\Omega$
such  that $\Lambda_1 f\big|_{\R\times\Gamma_0}=\Lambda f\big|_{\R\times\Gamma_0}$  for  all 
$f\in \R\times \Gamma_0$
where  $\Lambda_1,\Lambda$  are  two  DN  operators  for  $g_1$  and  $g$,  respectively.
Then,  as in  \S 2,   $g=g_1$  modulo  of the change  of  variables  (\ref{eq:1.7})  outside  of the  ergosphere.

Therefore,  without loss of  generality,  we  can assume  the  ergosphere  $\Delta=0$  for  metrics  $g,g_1$  is  the same  and  the  metrics  
$g,g_1$   are equal on $\Delta=0$.

Let    $\varphi^\pm$  be  a solution  of characteristic equation  of the form  (\ref{eq:3.2})  
with $[\hat g^{\rho\theta}]_{j,k=1}^2$  replaced  by  $[g_1^{\rho\theta}]_{j,k=1}^2$.

We assume,  as  in $(\ref{eq:4.2})$,   that 
\begin{equation}															\label{eq:4.7}
\varphi^\pm(\rho_1^\pm(x_0),\theta_1^\pm(x_0))=\varphi^\pm(\rho_1(0),\theta_1(0)),\ 
\theta_1(0)=\hat\theta_1.
\end{equation}

Make change  of variables  as  in  (\ref{eq:4.3})
\begin{equation}														\label{eq:4.8}
\sigma'=\varphi^+(\rho,\theta),\ \ \tau'=\varphi^-(\rho,\theta),
\end{equation}
where
\begin{equation}														\label{eq:4.9}
(\rho',\theta')\in \Pi',\ \ \Pi'=\{0\leq \rho'<\rho_0'(\theta'),\theta'\in \R^1\},
\end{equation}
$\rho=\rho_0'(\theta'),\theta' \in  \R$  is the black hole event horizon for the metric  $g_1$.
Note that  $L'u'=0$  has  the  form  $\frac{\partial^2 u'}{\partial\sigma'\partial\tau'}=0$   in  $(\sigma',\tau')$  coordinates.

As in  (\ref{eq:4.5}),  (\ref{eq:4.6}),   make  a change of variables                           
\begin{align}															\label{eq:4.10}
&
y_1'=\frac{\sigma'+\tau'}{2}=\frac{\varphi^+(\rho,\theta)+\varphi^-(\rho,\theta)}{2},
\\
\nonumber
&
y_2'=\frac{\sigma'-\tau'}{2}=\frac{\varphi^+-\varphi^-}{2},
\end{align}                                                     
\begin{equation}															\label{eq:4.11}
y_1'\big|_{\rho=0}=\theta,\ \ y_2'\big|_{\rho=0}=0.
\end{equation}
Note that  $L'u'=0$  has the form  $\big(\frac{\partial^2}{\partial y_1^2}-\frac{\partial^2}{\partial y_2^2}\big) u=0$,  $(y_1,y_2)\in \R_+^2$.

Let  $\Phi_1$  be  the map  (\ref{eq:4.10}),   $\Phi_1$  is  a homeomorphism  of  $\Pi'$ 
 onto  $\overline \R_+^2$.
Therefore,  $\Phi_0=\Phi^{-1}\Phi_1$  is a homeomorphism  of  $\Pi'$  onto $\Pi$.   Note  that   $\gamma_0'=\partial\Pi'$  
is  the black hole event horizon  for the metric $g_1$.   Since  the closure $\overline\Phi_0$  maps  $\overline\Pi'$  onto  $\overline\Pi$
we  have  that $\gamma_0'$  is mapped  onto $\gamma_0$.
Thus the event  horizon  $\gamma_0$  can be recovered  up to a change of variables.
Therefore the following theorem  holds:
\begin{theorem}															\label{theo:4.1}
Suppose  we have  two  wave equations  $Lu=o,\ L'u'=0$  such that  corresponding  DN  operators  
are equal  on  $\R\times\Gamma_0$.   Suppose  the ergosphere  $\Delta(x)=0$  of  $Lu=0$  is not  characteristic  for  any  $\Delta(x)=0$
and the trapped region  $O_\e\subset\Delta$  exists.  Then  the  ergosphere  $\Delta'(x)=0$  of  $L'u'=0$  is  also  
non-characteristic  for all  $\Delta'(x')=0$  and  has a trapped region $O_\e'$.  Moreover,  if  $\overline\Pi$  is the  black  hole of 
$Lu=0$  and  $\overline\Pi'$  is  the black  hole of  $L'u'=0$,  then  $\overline\Pi$  and $\overline\Pi'$  are  homeomorphic. 
\end{theorem}

{\bf Remark 4.1}   Since   the  kernel  of  $[g^{jk}]_{j,k=1}^2$  is  one-dimensional  for all $x$  belonging
to the ergosphere   $\Delta$,  there exists a vector  $e(x)\in \mbox{ker}\, [g^{jk}]_{j,k=1}^2$  smoothly  depending  on $x\in\Delta$.
Above we considered  two cases  when  $e(x)$  is normal to  $\Delta(x)=0$  for  all $x\in \Delta$  and  when  $e(x)$  is not  normal to $\Delta$  
for all $x\in\Delta$.   The last condition  is  equivalent  to  the  conditions
that the null-geodesics 
on $\Delta$  are not tangent to $\Delta$.

There is also a case  3)  when  $e(x)$  is normal  to  $\Delta=0$  only on some subset  of  $\Delta=0$.  Black  holes  in the case 
3)  were studied  in  [11]   and we shall explain  the ideas  behind  these studies.

The  underlying idea  in analyzing  the  black  holes  in two  space  dimensions   is the following:

Consider  the Hamiltonian
\begin{equation}														\label{eq:4.12}
H(x_1,x_2,\xi_0,\xi_1,\xi_2)=\sum_{j,k=0}^2 g^{jk}(x)\xi_j\xi_k.
\end{equation}
Let
\begin{align}														\label{eq:4.13}
&\frac{d x_k}{d s}=\frac{\partial H(x,\xi)}{\partial\xi_k},\ \ \frac{d\xi_k}{ds}=-\frac{\partial H}{\partial x_k},\ \ 0\leq k\leq 3,
\\  
\nonumber
& x_k(0)=y_k,\ \ \ \xi_k(0)=\eta_k
\end{align}
be  the equation  of  null-biocharacteristics.  Thus
\begin{equation}													\label{eq:4.14}
H\big(x_1(s),x_2(s),\xi_0(s),\xi_1(s),\xi_2(s)\big)=0\ \ \mbox{for  all} \ \  s.
\end{equation}
Since
$H$  is  independent  of  $x_0$  we have  that  $\frac{d\xi_0}{ds}=0$,  i.e.  $\xi_0(s)=\eta_0$  is  a constant.  We choose  
$\xi_0= 0$,  and  we shall call the null-bicharacteristics  with  $\xi_0=0$  the zero  energy  null-bicharacteristics.
The projection  of zero  energy  null-bicharacterstcs  on the  $(x_1,x_2)$-space  is called  the zero-energy  null-geodesics.  
Therefore  we have 
\begin{equation}													\label{eq:4.15}
\sum_{j,k=1}^2  g^{jk}(x(s))\xi_j(s)\xi_k(s)\equiv 0,\ \ \forall s.
\end{equation}
Equation (\ref{eq:4.15})  is  a quadratic equation  in $\xi_j(s),\ 1\leq j\leq 2$,  and therefore we have two families
of solutions
\begin{equation}															\label{eq:4.16}
\xi_j^\pm(s)=p_j^\pm(x(s)),\ \ x=(x_1,x_2),\ \ j=1,2.
\end{equation}
If  substitute  $\xi_j^\pm$   in (\ref{eq:4.13})  and  choose  $x_0$  as a parameter instead  of $s$  we obtain  two  $2\times 2$
system  of differential equations  in $(x_1,x_2)$:
\begin{equation}															\label{eq:4.17}
\frac{dx_j^\pm}{dx_0}=\frac{g^{j1}( x^\pm)p_1^\pm(x^\pm)+g^{j2}(x^\pm)p_2^\pm(x^\pm)}
{g^{01}( x^\pm)p_1^\pm(x^\pm)+g^{02}(x^\pm)p_2^\pm(x^\pm)},\ \ j=1,2.
\end{equation}
 Therefore  the solution of  $4\times 4$  system  of null-bicharacteristics  (\ref{eq:4.13})   is  reduced to the solution of $2\times 2$ 
 system  (\ref{eq:4.17}).   This reduction  substantially  simplifies  the study   of black  holes.  In the case  when  the boundary of
 the  ergoregion  is not characteristic  it was  done  in \S3.
 
 In the case  3)  black  holes  also  exist  and the boundary  of the  black hole  consists  of segments  of ``plus" or ``minus" 
 zero  energy  null-geodesics. In some cases  the boundary  of black  hole may have  corners  when  ``plus"
 null-geodesics  and  ``minus" null-geodesics  intersect.
 
 We will not  consider the inverse  problems  for  the black holes  in the case 3).

\section{The case  of  the axisymmetric metrics  in  three space  dimensions}
\init

Consider  the wave  equation  (\ref{eq:1.2})  when  $n=3$  and the metric  $g=\sum_{j,k=0}^3g_{jk} dx_j dx_k$  is  axisymmetric.

Let  $(\rho,\varphi,z)$  be  the  cylindrical  coordinates  $x=\rho\cos\varphi,\ y=\rho\sin\varphi,\ z=z.$ 
Since  metric  is axisymmetric  the Hamiltonian  of the spatial part of  (\ref{eq:1.2})  has  the form
\begin{equation}																\label{eq:5.1}
H=g^{\rho\rho}(\rho,z)\xi_\rho^2  +2g^{\rho z}(\rho,z)\xi_\rho\xi_z +g^{zz}\xi_z^2  +g^{\rho\varphi}\xi_\rho\xi_\varphi
+g^{z\varphi}\xi_z\xi_\varphi+g^{\varphi\varphi}\xi_\varphi^2,
\end{equation}
 where
the  coefficients of  (\ref{eq:5.1})  are independent of  $\varphi$.   Therefore  the dual variable  $\xi_\varphi$  is  constant.  
We choose  $\xi_\varphi=0,$ i.e.  instead of  (\ref{eq:5.1})  we have
\begin{equation}																\label{eq:5.2}
H_0=g^{\rho\rho}\xi_\rho^2+g^{\rho z}\xi_\rho \xi_z+g^{zz}(\rho,z)\xi_z^2.
\end{equation}
We shall call (\ref{eq:5.2})  the restricted  axisymmetric Hamiltonian.  The wave equation    
corresponding  to the  restricted  axisymmetric  Hamiltonian  
is obtained  from (\ref{eq:1.2})  by  dropping all derivatives  in $\varphi$.

Factoring  $H_0=0$   we get
\begin{equation}																\label{eq:5.3}
\xi_\rho=\frac{
-g^{\rho z}\pm 
\sqrt
{
(g^{\rho z})^2-g^{\rho\rho} g^{z z}
}
}
{g^{\rho\rho}
}
\xi_z.
\end{equation}
Note that  
\begin{equation}																\label{eq:5.4}
\Delta_1(\rho,z)=g^{\rho\rho}g^{zz}-(g^{\rho z})^2\leq 0
\end{equation}
is the  ergoregion    for  the  restricted  axisymmetric  equation.   
As  in \S3  there exists  a smooth vector  $e_1(\rho,z)$  on  $\Delta_1(\rho,z)=0$   such
that
\begin{equation}																\label{eq:5.5}
\left[\begin{matrix}
 g^{\rho\rho}   &   g^{\rho z}\\
 g^{\rho z} &   g^{zz}
\end{matrix} \right]
e_1(\rho,z)=0 \ \ \mbox{for  all}\ \ (\rho,z)\in \Delta_1=0.
\end{equation}
As  in \S3  there are  three choices:
\begin{align}																	\label{eq:5.6}
&
a)\ \   e_1(\rho,z)=C\Big(\frac{\partial \Delta_1}{\partial\rho},\frac{\partial\Delta_1}{\partial z}\Big) \ \ \mbox{for all}\ \ (\rho,z)\in\{ \Delta_1=0\},
\\
\nonumber
&
b)\ \  e_1(\rho,z)\neq C\Big(\frac{\partial \Delta_1}{\partial\rho},\frac{\partial\Delta_1}{\partial z}\Big) \ \ \mbox{for  all}\ \ 
(\rho,z)\in\{ \Delta_1=0\} \ \ \mbox{and}
\\
\nonumber
&
c)\ \   e_1(\rho,z)=C\Big(\frac{\partial \Delta_1}{\partial\rho},\frac{\partial\Delta_1}{\partial z}\Big)\ \ \
\mbox{only for a subset  of }\ \ (\rho,z)\in \{\Delta_1=0\}.
\end{align}
In the case  a)  $\Delta_1(\rho,z)=0$  is  a characteristic curve and  we have a Schwarzschield type  black hole.
In the case  b)  there  exists   a black  hole  inside  $\{\Delta_1(\rho,z)<0\}$.
The case  c)   will not be considered  in this paper.

Consider  the case of  Kerr metric  -  the most celebrated   example  of the axisymmetric  metric  (cf. [18],  [23]).  The
restricted Kerr Hamiltonian 
   is  (cf.  (\ref{eq:5.2})):
\begin{equation}																	\label{eq:5.7}
H_0(\rho,z,\xi_\rho,\xi_z)=(-\xi_\rho^2-\xi_z^2)+K(b_\rho\xi_\rho+b_z\xi_z)^2,
\end{equation}
where
\begin{equation}																\label{eq:5.8}
K=\frac{2mr^3}{r^4+a^2z^2},\ \ b_\rho=\frac{\rho^2}{r^2+a^2},\ \ b_z=\frac{z}{r},
\end{equation}
where  $r$  is defined  by the relation
\begin{equation}																\label{eq:5.9}
\frac{\rho^2}{r^2+a^2}+\frac{z^2}{r^2}=1.
\end{equation}
We  have
\begin{equation}																	\label{eq:5.10}
g^{\rho\rho}=-1+Kb_\rho^2,
\ \
g^{zz}=-1+Kb_z^2,
\ \
g^{\rho z}=2Kb_\rho b_z,
\end{equation}
where  $g^{\rho\rho}, g^{\rho z},g^{zz}$  are   the same  as  in  (\ref{eq:5.2}).

Thus 
\begin{equation}  																\label{eq:5.11}
\hat\Delta_1(\rho,z)=\det \left[\begin{matrix}
 g^{\rho\rho}   &   g^{\rho z}\\
 g^{\rho z} &   g^{zz}
\end{matrix} \right]=Kb_\rho^2+Kb_z^2-1.
\end{equation}

It was shown  in  [10]  that
\begin{equation}																\label{eq:5.12}
\hat\Delta_1(\rho,z)=C_\pm(r-r_\pm),
\end{equation}
where  
$$
C_\pm\neq 0,\  r_\pm=m\pm \sqrt{m^2-a^2}.
$$

Thus  there  are  two  ergospheres:  the outer  ergosphere  $\Delta_1^+(\rho,z)=C_0^+(r-r_+)=0 $  and the inner  ergosphere  
$\Delta_1^-=C_0^-(r-r_-)=0$   for the restricted  Kerr metric.
It is  known  (cf. [23])  that  
\begin{equation}																\label{eq:5.13}
\frac{\rho^2}{r_\pm^2+a^2}+\frac{z^2}{r_\pm^2}=1
\end{equation}
are characteristic  surfaces  for  the Kerr  equation.

Note that
\begin{align}																\label{eq:5.14}
&r-r_+=C_1^+\Big(\frac{\rho^2}{r_+^2+a^2}+\frac{z^2}{r_+^2}-1\Big),
\\																		
\label{eq:5.15}
&r-r_-=C_1^-\Big(\frac{\rho^2}{r_-^2+a^2}+\frac{z^2}{r_-^2}-1\Big),
\end{align}
where  $C_1^+,C_1^-\neq 0$.

Therefore
$\hat\Delta_1=\hat\Delta_1^+\hat\Delta_1^-$.
\begin{equation}															\label{eq:5.16}
\hat\Delta_1^\pm(\rho,z)=C_3^\pm\Big(\frac{\rho^2}{r_\pm^2+a^2}+\frac{z^2}{r_\pm^2}-1\Big).
\end{equation}
Denote
\begin{equation}															\label{eq:5.17}
S^\pm=\frac{\rho^2}{r_\pm^2+a^2}+\frac{z^2}{r_\pm^2}-1.
\end{equation}
Then
\begin{align}																\label{eq:5.18}
&\frac{\partial}{\partial\rho}\hat\Delta_1^\pm=\Big(\frac{\partial}{\partial\rho}C_3^\pm\Big)\Big(\frac{\rho^2}{r_\pm^2+a^2}+
\frac{z^2}{r_\pm^2}-1\Big)
+C_3^\pm\frac{\partial}{\partial \rho}\Big(\frac{\rho^2}{r_\pm^2+a^2}+\frac{z^2}{r_\pm^2}-1\Big),
\\
																		\label{eq:5.19}
&\frac{\partial}{\partial z}\hat\Delta_1^\pm=\Big(\frac{\partial}{\partial z}C_3^\pm\Big)\Big(\frac{\rho^2}{r_\pm^2+a^2}+\frac{z^2}{r_\pm^2}-1\Big)
+C_3^\pm\frac{\partial}{\partial z}\Big(\frac{\rho^2}{r_\pm^2+a^2}+\frac{z^2}{r_\pm^2}-1\Big).										
\end{align}								
When  																		
$S^\pm=\frac{\rho^2}{r_\pm^2+a^2}+\frac{z^2}{r_\pm^2}-1=0$ we got  that the  first terms  in  (\ref{eq:5.18}),  (\ref{eq:5.19})  vanish.

Therefore
\begin{align}
\nonumber
&g^{\rho\rho}\Big(\frac{\partial}{\partial\rho}\hat\Delta_1^\pm\Big)^2
+2 g^{\rho\theta}\frac{\partial \hat\Delta_1^\pm}{\partial\rho}\frac{\partial\hat\Delta_1^\pm}{\partial z} + g^{\theta\theta}
\Big(\frac{\partial\hat\Delta_1^\pm}{\partial z}\Big)^2
\\
\nonumber
= &g^{\rho\rho}\Big(S_\rho^\pm\Big)^2+2g^{\rho\theta}\frac{\partial S^\pm}{\partial\rho}\frac{\partial S^\pm}{\partial z}
+ g^{\theta\theta}\Big(\frac{\partial S^\pm}{\partial z}\Big)^2=0
\end{align}
since   $S^\pm=0$  and  $S^\pm$  are   characteristic  surfaces.

Therefore  the following  theorem  holds:
\begin{theorem}															\label{theo:5.1}
We have that
  $\hat\Delta_1^+=0$  is an outer  ergosphere  and an outer black hole horizon  and  $\hat\Delta_1^-=0$
is        an inner  ergosphere  and  an inner black hole  horizon  for the restricted Kerr equation'
\end{theorem}


\end{document}